\documentclass[12pt]{iopart}

\newtheorem{theorem}{Theorem}{}
\newtheorem{corollary}{Corollary}{}
\newtheorem{remark}{Remark}{}
\newtheorem{lemma}{Lemma}{}
\newtheorem{definition}{Definition}
\newtheorem{prop}{Proposition}

\usepackage{graphicx}

\begin{document}

\title[]{A Necessary and Sufficient Condition for Complete phase synchronization of high-dimensional nonidentical Kuramoto   oscillators}

\author{Yushi Shi, Ting Li \& Jiandong Zhu}

\address{School of Mathematical
Sciences, Nanjing Normal University, Nanjing, China}
\ead{zhujiandong@njnu.edu.cn}
\vspace{10pt}
\begin{indented}
\item[]August 2022
\end{indented}

\begin{abstract}
For original Kuramoto models with nonidentical oscillators, it is impossible to realize complete phase synchronization. However, this paper reveals that complete phase synchronization  can be achieved for a large class of high-dimensional Kuramoto models with nonidentical oscillators. Under the topology of strongly connected digraphs, a necessary and sufficient condition for complete phase  synchronization is obtained. Finally, some simulations are provided to validate the obtained theoretic results.   
\end{abstract}

%
%
%
%
%

\section{Introduction}
\ \ \ \ \ \
The collective behavior of multi-agent system is a hot topic in complex system research. A considerable number of individuals have their own dynamic rules, individuals interact with each other through local information exchange. These interacting individuals are also called coupled oscillators. Synchronization is a fundamental problem in the study of collective behavior, which widely exists in nature, society, physics, biology and other systems\cite{couzin2018synchronization,shahal2020synchronization,Dorfler,Grzybowski,cumin2007,Simone-cta}. The famous Kuramoto model is a nonlinear system model used to describe the synchronization behavior of a large number of coupled oscillators, it was first proposed by Yoshiki Kuramoto in 1975 \cite{kuramoto1975}.The dynamic equation of the original Kuramoto model with $m$ oscillators is described as:
\begin{equation}
\label{eq1}
\dot \theta_{i}=\omega_{i}+k\sum\limits_{j=1}^ma_{ij}\sin(\theta_{j}-\theta_{i}),  \ \ i=1,2,\ldots,m,
\end{equation}
where $\theta_{i}$ is the phase of the oscillator $i$, $\omega_i$ is the
natural frequency, $k>0$ is the coupling gain, and $A=(a_{ij})\in{\mathbf{R}^{m\times{m}}}$ is the weighted adjacency matrix of the information communication graph $\mathcal{G}$. Olfati-Saber \cite{Saber} and Lohe \cite{lohe2009,lohe2010} independently generalized the traditional Kuramoto model to the general high-dimensional form named high-dimensional Kuramoto model or Lohe model.
\par
The Kuramoto model with $r_{i}$ as the state of the $i$th oscillator in the $n$-dimensional real linear space $\mathbf{R}^{n}$ is described as follows:
\begin{equation}
\label{eq01}
\dot{r}_{i}=\Omega_{i}r_{i}+k\sum\limits_{j=1}^{m}a_{ij}\big(r_{j}-\frac{r_{i}^{\mathrm{T}}r_{j}}{r_{i}^{\mathrm{T}}r_{i}}r_{i}\big), \ \  i=1,2,\ldots,m,
\end{equation}
where $\Omega_{i}$ is a real $n\times n$ skew-symmetric matrix. When $\Omega_{i}=\Omega$ for each $i=1,2,\ldots,m,$ (\ref{eq01}) is called the high-dimensional Kuramoto model with identical oscillators, otherwise, (\ref{eq01}) is called the high-dimensional Kuramoto model with nonidentical oscillators or high-dimensional nonidentical Kuramoto   oscillators. In our earlier work \cite{zhu2013}, we have been proved that, under the condition of $\Omega_{i}$ is a skew-symmetric matrix, the value of $\|r_{i}(t)\|$ is a constant. Without loss of generality, we give the following high-dimensional Kuramoto model restricted to the unit sphere:
 \begin{equation}
 \label{eq001}
 \dot{r}_{i}=\Omega_{i}r_{i}+k\sum\limits_{j=1}^{m}a_{ij}\big(r_{j}-{(r_{i}^{\mathrm{T}}r_{j})}r_{i}\big), \ \ \|r_i\|=1,\ \  i=1,2,\ldots,m.
 \end{equation}

\par
With respect to (\ref{eq01}) or (\ref{eq001}) with identical oscillators and its variations, there have been many theoretical analyses and derivations on the synchronization problem  \cite{Ha2016,chi2014,zhang2018,Markdahl-Automatica-2018,zhu2014high,HA2021132781,MARKDAHL2020108736}.
In \cite{zhangarxiv}, exponential synchronization was proved under some initial state constraints. In our recent paper \cite{peng2021exact}, the exact exponential synchronization rate has been revealed. For nonidentical oscillators cases, the synchronization problem of Kuramoto model on the two-dimensional plane has been a large number research results \cite{ha2015practical,PhysRevE022302,wu8338121,gil2021optimally,clusella2022kuramoto}. However, the research results of synchronization problem in high dimensional case are relatively few. Some literatures only study the information connection graph in the case of complete graph, and get the practical synchronization results under some initial conditions instead of complete synchronization\cite{choi2014,choi2016}. In \cite{Markdahl2021CP}, Markdahl et al. first revealed that a specific class of high-dimensional nonidentical Kuramoto oscillators can achieve complete phase synchronization. But the graph considered in \cite{Markdahl2021CP} is undirected, and the condition on the frequency matrices is only a sufficient condition. In \cite{zhang2021synchronization}, for proportional nonidentical Kuramoto oscillators interconnected by a digraph, it was revealed that the complete phase synchronization can occur, and a rigorous theoretical analysis was provided. However, the complete phase synchronization problem of the general high-dimensional Kuramoto model with nonidentical oscillators is still an open problem.\par
In this paper, for the high-dimensional Kuramoto model with nonidentical oscillators interconnected by a strongly connected digraph, a necessary and sufficient condition for  complete phase synchronization is proposed. Some numerical examples are given to validate the obtained theoretical results.
\par
The rest of this paper is organized as follows. Section 2 gives some  preliminaries and the problem statement. Section 3 includes our main result. Section 4 gives some simulations. Finally, section 5 is devoted to a summary.

\section{Preliminaries and problem statement}

Consider the autonomous system

\begin{equation}
\label{eqn001}
\dot x=f(x),
\end{equation}
where $f:D\rightarrow \mathbf{R}^n$ is a locally Lipschitz
map from a domain $D\subset \mathbf{R}^n$ into $\mathbf{R}^n$.

\begin{definition}\cite{Khalil}
Let $x(t)$ be a solution of (\ref{eqn001}). A point \(\eta \) is said to be a positive limit point of $x(t)$ if there is a sequence $\{t_n\}$, with $t_n\rightarrow +\infty$ as $n\rightarrow +\infty$, such that $x(t_n)\rightarrow \eta$ as $n\rightarrow +\infty$. The set of all positive limit points of $x(t)$ is called the positive limit set of $x(t)$. A set $M$ is said to be an invariant set with respect to (\ref{eqn001}) if
$$
x(0)\in M \Rightarrow x(t)\in M,\ \ \forall \ t\in R.
$$
A set $M$ is said to be a positively invariant set with respect to (\ref{eqn001}) if
$$
x(0)\in M \Rightarrow x(t)\in M,\ \ \forall \ t\geq 0.
$$
\end{definition}

\begin{lemma}
\label{lem1} (Lemma 4.1 of \cite{Khalil})
If a solution $x(t)$ of (\ref{eqn001}) is bounded and belongs to domain $D$ for $t\geq 0$, then its positive limit set $L^+$ is a nonempty, compact, invariant set. Moreover, $x(t)$ approaches $L^+$ as $t\rightarrow +\infty$.
\end{lemma}
\begin{prop}
\label{prop1}
(LaSalle's Theorem, Theorem 4.4 of \cite{Khalil})
Let $C\subset D$ be a compact set that is positively invariant with respect to (\ref{eqn001}). Let $V:D\rightarrow \mathbf{R}$ be a continuously differentiable function such that $\dot V(x)\leq 0$ in $C$. Let $E$ be the set of all points in $C$ where $\dot V(x)=0$. Let $M$ be the largest invariant set in $E$. Then every solution starting in $C$ approaches $M$ as $t\rightarrow +\infty$.
\end{prop}

The high-dimensional Kuramoto model can be regarded as a special nonlinear multi-agent system. Since the interconnecting graph is the fundamental element of a multi-agent system, we first introduce some relevant knowledge about algebraic graph theory.\par
Let $\mathcal{G}=(\mathcal{V},\mathcal{E},A)$ be a weighted digraph, where  $\mathcal{V}=\{1,2,\ldots,m\}$ is the set of agents,  $\mathcal{E}\subset{\mathcal{V}\times \mathcal{V}}$ is the set of edges, and  $A=(a_{ij})$ is the weighted adjacent matrix, satisfying
\begin{equation*}
a_{ij}\left\{\begin{array}{ll}
>0,& \quad (j,i)\in{\mathcal{E}},\\
=0 ,&  \quad \mathrm{otherwise}.
\end{array}\right.
\end{equation*}
A directed edge $(i,j)$ means that agent $j$ can receive the state information from agent $i$.
If any two agents $x,y\in{\mathcal{V}}$ can get information about each other through the directed edges, digraph $\mathcal{G}$ is said to be strongly connected.
The Laplacian matrix $L=(l_{ij})$ of $\mathcal{G}$ is defined by
\begin{equation*}
l_{ij}=\left\{\begin{array}{ll}
-a_{ij},& \quad i\neq j,\\
\sum\limits_{k\neq i}a_{ik}, & \quad i=j.
\end{array}\right.
\end{equation*}
\begin{lemma}(Perron-Frobenius Theorem, Theorem 8.8.1 of \cite{Godsil2001})
\label{lemma2000}
Suppose that $A$ is a real nonnegative square matrix whose underlying directed graph is strongly connected. Then the spectral radius $\rho(A)$ is a simple eigenvalue of $A$. If $x$ is an eigenvector for  $\rho(A)$, then no entries of $x$ are zero, and all have the same sign.
\end{lemma}

\begin{lemma}
\label{lemma201} (Corollary 3 of \cite{scutari2008}) Let
$\mathcal{G}=(\mathcal{V},\mathcal{E},A)$ be a digraph with
Laplacian matrix $L$. If $\mathcal{G}$ is strongly connected, then
$L$ has a simple zero eigenvalue and a positive left-eigenvector
associated to the zero eigenvalue.
\end{lemma}





\begin{definition}\cite{choi2016,ha2019emergent} 
\label{def1}
Let $r(t)=(r_{1}^{\mathrm{T}},\ldots,r_{m}^{\mathrm{T}})^{\mathrm{T}}$ be a solution of system (\ref{eq001}). It is said that the {\it complete phase synchronization} is achieved if
\begin{equation}
\label{eqn005}
\lim_{t\rightarrow +\infty}(r_i(t)-r_j(t))=0,\ \ \forall \ 1\leq i,\  j\leq  m.
\end{equation}
\end{definition}

In particular, when $n=2$, we let$$
\Omega_i=\left[\begin{array}{cc}
               0 & -\omega_i\\
               \omega_i &0
               \end{array}\right],\ \ r_i=\left[\begin{array}{c}
               \cos\theta_i\\
                \sin\theta_i\\
                             \end{array}\right].
$$
 Then the dynamics of $\theta_i$ is just the original Kuramoto model. When $n=2$,  it is impossible to achieve complete phase synchronization for nonidentical oscillators.

A natural question is whether the high-dimensional nonidentical Kuramoto oscillators  can achieve complete phase synchronization when \(n\geq 3\). Actually, this question has been partially answered in \cite{Markdahl2021CP} and \cite{zhang2021synchronization}  for some special class of
high-dimensional nonidentical oscillators. The main task of this paper is to explore  necessary and sufficient conditions for solving this problem.


\section{Main Results }
\ \ \ \ \ \ We first investigate some necessary conditions for complete phase synchronization.
\begin{theorem}
        \label{thm1}
                Consider the high-dimensional Kuramoto model with nonidentical
                oscillators limited to the unit sphere and described by (\ref{eq001}). If there is a solution
$$
r(t)=(r_1^{\mathrm T}(t),\ r_2^{\mathrm T}(t),\ldots,  r_m^{\mathrm T}(t))^{\mathrm T}\in \mathbf{R}^{mn}
$$
of (\ref{eq001}) exhibits complete phase synchronization, then $\dim(W)>0$ (or equivalently \(W\ne \{0\}\)), where
   \begin{equation}
\label{eqn006}
W=\bigcap_{l=1}^{n}\{\xi\ |\  \Omega_1^l\xi=\Omega_2^l\xi=\cdots=\Omega_m^l\xi\}=\bigcap_{l=1}^{n}\bigcap_{i=2}^m\ker(\Omega_i^l-\Omega_1^l).
\end{equation}
\end{theorem}
{\bf Proof}. Since $\|r_i(t)\|=1$ for all $i=1,2,\ldots,m$, the solution $r(t)$ is bounded. Thus, by Lemma \ref{lem1},  the positive limit set $L^+$ of $r(t)$ is a nonempty, compact, invariant set. Since \(r(t) \) exhibits complete phase synchronization on the unit sphere, every point $\eta\in L^+$ has the form
$$
\eta=(\bar\eta^{\mathrm T}, \bar\eta^{\mathrm T},\ldots,\bar\eta^{\mathrm T})^{\mathrm T}\in \mathbf{R}^{mn}, \ \ \ \|\bar\eta\|=1 .
$$
Considering the invariance of $L^+$, we write a solution of $(\ref{eq001})$ in $L^+$ into the form
$$
\eta(t)=(\bar\eta^{\mathrm T}(t),\ \bar\eta^{\mathrm T}(t),\ldots,\bar\eta^{\mathrm T}(t))^{\mathrm T}\in \mathbf{R}^{mn},\ \ \|\bar\eta(t)\|=1.
$$
Substituting the solution $\eta(t)$ into $(\ref{eq001})$ yields $\dot {\bar\eta}(t)=\Omega_i \bar\eta(t)$ for every $i=1,2,\ldots,m$. Taking the $l$-th order derivative with respect to $t$, we have ${\bar\eta}^{(l)}(t)=\Omega_i^l{\bar\eta}(t)$.  which implies $\Omega_1^l{\bar\eta}(t)=\Omega_2^l{\bar\eta}(t)=\cdots=\Omega_m^l{\bar\eta}(t)$. Further considering $\|\bar\eta(t)\|=1$, we conclude $W\ne \{0\}$.

\begin{corollary}
\label{cor1}
If $W=\{0\}$, then the high-dimensional nonidentical Kuramoto oscillators described by (\ref{eq001}) cannot achieve complete phase synchronization.
\end{corollary}

\begin{corollary}
\label{cor2}
If there exist $\Omega_i$ and $\Omega_j$ such that $\det(\Omega_i-\Omega_j)\ne 0$, then the high-dimensional nonidentical Kuramoto oscillators described by (\ref{eq001}) cannot achieve complete phase synchronization.
\end{corollary}

\begin{corollary}
\label{cor3}
If $n=2$, then the nonidentical Kuramoto oscillators described by (\ref{eq001}) cannot achieve complete phase synchronization.
\end{corollary}
{\bf Proof}. Since $n=2$, without loss of generality, we assume that
$$
\Omega_1=\left[\begin{array}{cc}
               0 & -\omega_1\\
               \omega_1 &0
               \end{array}\right],\ \ \Omega_2=\left[\begin{array}{cc}
               0 & -\omega_2\\
               \omega_2 &0
               \end{array}\right],\ \ \omega_1\ne\omega_2.
$$
Obviously, $\det(\Omega_1-\Omega_2)\ne 0$. So the proof is complete by Corollary \ref{cor2}.

\begin{lemma}
\label{lem5}
The linear subspace $W$ shown in (\ref{eqn006}) of Theorem \ref{thm1} has the alternative expressions as follows:
\begin{equation}
\label{eqn07}
W=\bigcap_{s=0}^{n-1}\bigcap_{i=2}^m\ker\Big((\Omega_i-\Omega_1)\Omega_1^s\Big)=\bigcap_{s=0}^{\infty}\bigcap_{i=2}^m\ker\Big((\Omega_i-\Omega_1)\Omega_1^s\Big).
\end{equation}
\end{lemma}
{\bf Proof}. For any $\xi\in W$, we have $\Omega_1^l\xi=\Omega_i^l\xi$ for every $l=1,2,\ldots,n$ and
$i=1,2,\ldots,m$. Thus, for any $s=0,1,2,\ldots,n-1$, we have
\begin{equation}
\!\!\!\!\!\!\!\!\!(\Omega_i-\Omega_1)\Omega_1^s\xi=\Omega_i\Omega_1^s\xi-\Omega_1^{s+1}\xi =\Omega_i\Omega_i^s\xi-\Omega_1^{s+1}\xi=\Omega_i^{s+1}\xi-\Omega_1^{s+1}\xi=0.
\end{equation}
Conversely, assume $(\Omega_i-\Omega_1)\Omega_1^s\xi=0$ for any $s=0,1,2,\ldots,n-1$. Then, when $s=0$, it follows that $(\Omega_i-\Omega_1)\xi=0$, that is, $\Omega_i\xi=\Omega_1\xi$. When $s=1$, we have $(\Omega_i-\Omega_1)\Omega_1\xi=0$,
which implies that
$$
\Omega_1^2\xi=\Omega_i\Omega_1\xi=\Omega_i\Omega_i\xi=\Omega_i^2\xi.
$$
Letting $s=2$, we obtain $(\Omega_i-\Omega_1)\Omega_1^2\xi=0$, which results in
$$
\Omega_1^3\xi=\Omega_i\Omega_1^2\xi=\Omega_i\Omega_i^2\xi=\Omega_i^3\xi.
$$
With the same argument, one can conclude that
$\Omega_1^l\xi=\Omega_i^l\xi$ for every $l=1,2,\ldots,n$.
Therefore, the first equality of $(\ref{eqn07})$ is proved, and the second equality of $(\ref{eqn07})$ is subsequently obtained by  Hamilton-Cayley Theorem.
\begin{remark}
\label{rek1}
We have known that, when $n=2$,  the high-dimensional Kuramoto oscillators  described by (\ref{eq001}) reduce to the original Kuramoto oscillators (\ref{eq1}). Therefore, Corollary \ref{cor3} means that, if complete phase synchronization of the original Kuramoto model (\ref{eq1}) is achieved, then the oscillators have a common natural frequency, which is just a well-known result on the original Kuramoto model \cite{Lunze2011}. When $n=2$, the condition $\dim (W)>0$ is equivalent to $\Omega_i=\Omega_{j}$ for any $i$ and $j.$ In this case, the system is reduced to the identical Kuramoto model.
 \end{remark}

\begin{remark}
\label{rek021}
If the high-dimensional Kuramoto oscillators are identical, i.e. $\Omega_i=\Omega$ for all $i=1,2,\ldots, m$, then $W=\mathbf{R}^n$.
If the oscillators are proportional and nonidentical, i.e.  $\Omega_i=l_i\Omega$ for all $i=1,2,\ldots, m$ and there exist $i$ and $j$ such that $l_i\ne l_j$, then $W=\ker(\Omega)$. If $\Omega_i\Omega_1=\Omega_1\Omega_i$ for $i=2,\ldots, m$, it follows from Lemma \ref{lem5} that $W=\bigcap_{i=2}^m\ker(\Omega_i-\Omega_1)$.
 \end{remark}

In the above analysis, we have obtained some necessary conditions for  complete phase synchronization. In the following, we explore sufficient conditions and the attracting region for complete phase synchronization.

\begin{theorem}
        \label{thm2}
                Consider the high-dimensional Kuramoto model with nonidentical
                oscillators described by (\ref{eq001}). Assume that each $\Omega_i$ is skew-symmetric,  the linear subspace $W$ described by (\ref{eqn006}) satisfies  $\dim (W)>0$, and the interconnection graph $\mathcal{G}$ is a strongly connected digraph with the adjacency matrix $A=(a_{ij})$ . If the initial states satisfy $p^{\mathrm{T}}r_{i}(0)>0$ for some $p\in W$ and  every $i=1,2,\ldots,m$, then the solution $r(t) $  exhibits complete phase synchronization, and each \(r_{i}(t)\) approaches $W$ as $t\rightarrow +\infty$.
\end{theorem}
{\bf Proof}.
        We apply the LaSalle's Theorem to the autonomous system as follows:
\begin{eqnarray}
\label{eqn007}
 \dot{r}_{i}&=&\Omega_{i}r_{i}+k\sum\limits_{j=1}^{m}a_{ij}\big(r_{j}-{(r_{i}^{\mathrm{T}}r_{j})}r_{i}\big),\\ \label{eqn008}
\dot \eta&=& \Omega_1\eta,
\end{eqnarray}
where $r_i,\eta\in \mathbf{R}^n$.   For any given  $\varepsilon>0$, we construct a compact set
$$
C_\varepsilon\!=\!\!\{\!(r,\eta)|\ r\!=\!(r_1^{\mathrm T},\ldots,r_m^{\mathrm T})^{\mathrm T},\ \|r_i\|\!=\!1,\ \|\eta\|\!=\!1,\ \eta^{\mathrm{T}}r_{i}\!\geq \varepsilon,\ \eta\in W,\ i\!=\!1,2,\ldots,m\}.
$$

First, we show that $W$ is invariant with respect to (\ref{eqn008}). As a matter of fact, if     $\eta(0)\in W$,
then $(\Omega_i-\Omega_1)\Omega_1^l\eta(0)=0$ for any $1\leq i\leq m$ and $l\geq 0$ by Lemma \ref{lem5}.
Thus,
\begin{equation}
\!\!\!\!\!(\Omega_i-\Omega_1)\Omega_1^l\eta(t)=(\Omega_i-\Omega_1)\Omega_1^l\mathbf{e}^{\Omega_1t}\eta(0)=\sum_{i=0}^{\infty}\frac{1}{i!}(\Omega_i-\Omega_1)\Omega_1^{l+i}\eta(0)= 0,
\end{equation}
which implies that $\eta(t)\in W$.

Next, we prove that $C_{\varepsilon}$ is positively invariant. Let
$$
h_{i}(t)= \eta^\mathrm{T}(t)r_{i}(t),\ \  h(t)\!=\!\!\min\limits_{1\leq i\leq m}\!\!h_{i}(t), ~~v_{t}\!=\!\!\{i|h_{i}(t)\!=\!h(t)\}.
$$
Assuming $(r(t),\eta(t))\in C_{\varepsilon}$, which implies each $h_{i}(t)\geq \varepsilon$ and $\eta(t)\in W,$ we have
        \begin{eqnarray}
        \label{eqn010}
       \!\! \!\!\!\!\!\dot{h}_{i}(t)&=&\eta^\mathrm{T}(t)(\Omega_i-\Omega_1)r_{i}(t) +k\sum\limits_{j=1}^{m}a_{ij}(h_{j}(t)-(r_{i}^{\mathrm{T}}(t)r_{j}(t))h_{i}(t)) \nonumber \\
                &=&k\sum\limits_{j=1}^{m}a_{ij}\big(h_{j}(t)-(r_{i}^{\mathrm{T}}(t)r_{j}(t))h_{i}(t)\big)
                \end{eqnarray}
and
        \begin{eqnarray*}
                D^{+}h(t)&=&\min\limits_{i\in v_{t}}k\sum\limits_{j=1}^{m}a_{ij}(h_{j}(t)-(r_{i}^{\mathrm{T}}(t)r_{j}(t))h_{i}(t)). \\
                &\\
                &\geq& \min\limits_{i\in{v_{t}}}k\sum\limits_{j=1}^{m}a_{ij}(h_{j}(t)-h_{i}(t))\\
                &\geq&  \min\limits_{i\in{v_{t}}}k\sum\limits_{j=1}^{m}a_{ij}(h_{j}(t)-h(t))\\
                &\geq&0,
        \end{eqnarray*}
        where, for the operation rule of the Dini derivative, please refer to Lemma 2.2 of \cite{lin2007state}. So,  $h(t)$ is a nondecreasing function, and  subsequently all  $\eta^{\mathrm{T}}(0)r_{i}(0)\!\geq \varepsilon$ imply $\eta^{\mathrm{T}}(t)r_{i}(t)\!\geq \varepsilon$ for all $t\geq 0$.  Therefore, $C_{\varepsilon}$ is a positively invariant compact set of the system described by (\ref{eqn007}) and (\ref{eqn008}).

        Then, we define a function $V(r,\eta)$ in the positively invariant compact  set $C_{\varepsilon}$. Write (\ref{eqn010}) as
        \begin{eqnarray}
        \label{eqna3}
        \dot{h}_{i}&=&k\sum\limits_{j=1}^{m}a_{ij}(h_{j}-h_{i})+k\sum\limits_{j=1}^{m}a_{ij}(1-(r_{i}^{\mathrm{T}}r_{j}))h_{i},       \end{eqnarray}
        which can be rewritten in the compact form
        \begin{eqnarray}
        \label{eqna4}
        \dot{h}=-kLh+kF(r,h)
        \end{eqnarray}
        with $h=(h_{1},h_{2},\ldots,h_{m})^{\mathrm{T}}$, $L$ being the Laplacian matrix of $\mathcal{G}$ and
        $$
        F(r,h)=\Big(
        \sum\limits_{j=1}^{m}a_{1j}(1-r_{j}^{\mathrm{T}}r_{1})h_{1},
        ~\ldots~,\sum\limits_{j=1}^{m}a_{mj}(1-r_{j}^{\mathrm{T}}r_{m})h_{m}
        \Big)\!^\mathrm{T}.
        $$
        Since digraph $\mathcal{G}$ is strongly connected, by Lemma \ref{lemma201}, the Laplacian matrix $L$ has a positive left eigenvector associated with the zero eigenvalue, i.e. there exists vector $\beta=(\beta_{1},\beta_{2},\ldots,\beta_{m})^{\mathrm{T}}\in{\mathbf{R}^{m}}$ such that $\beta^{\mathrm{T}}L=0$ with all $\beta_{i}>0, i=1,2,\ldots,m$. Now, let us define $$V(r,\eta)=-\sum_{i=1}^m\beta_i\eta^\mathrm{T}r_{i}=-\beta^{\mathrm{T}}h.$$ It is straightforward to  validate that
        \begin{eqnarray}
        \label{eqna5}
        \dot{V}(r,\eta)&=&-\beta^{\mathrm{T}}\dot{h} \nonumber \\
        &=&-k\beta^{\mathrm{T}}Lh-k\beta^{\mathrm{T}}F(r,h) \nonumber \\
        &=&-k\sum\limits_{i,j=1}^{m}a_{ij}\beta_{i}(1-r_{i}^{\mathrm{T}}r_{j})h_{i} \nonumber \\
        &\leq& 0,\ \ \ \ \ \ \ \ \ \ \ \ \forall\  (r,\eta)\in C_\varepsilon.
        \end{eqnarray}
        In the following, we assume that the equality in (\ref{eqna5}) holds. Then
        \begin{eqnarray}
        \label{eqna6}
        a_{ij}\beta_{i}(1-r_{i}^{\mathrm{T}}r_{j})h_{i}=0\qquad (i,j=1,2,\ldots,m).
        \end{eqnarray}
        Since $\beta_{i}>0$ and $h_{i}>0$, it follows from (\ref{eqna6}) that $a_{ij}>0$ implies $r_{i}=r_{j}$. Therefore, by the strong connectedness of graph $\mathcal{G}$, we have $$E=\{(r,\eta)\in C_\varepsilon\ |\ \dot{V}(r,\eta)=0\}=\{(r,\eta)\in C_\varepsilon\ |\ r_{1}=r_{2}=\cdots=r_{m}\}.
$$

After that, let us investigate the largest invariant set $M$ in $E$. For any $( \tilde r, \tilde \eta)\in M$, it follows from the invariance of \(M\) that the  solution $(\tilde r(t), \tilde \eta(t))$ of (\ref{eqn007}) and (\ref{eqn008}) starting from the initial state $( \tilde r, \tilde \eta)$ is also in \(M\), which implies that $\tilde r_i(t)=\tilde r_1(t)$ and $\dot{\tilde r}_1(t)=\Omega_i {\tilde r}_{1}(t) $ for each $i=1,2,\ldots,m$. So, taking the $l$-th derivative and then letting \(t=0\), we have $\Omega_i^l{\tilde r_{1}}=\Omega_1^l{\tilde r_{1}}$. Therefore, ${\tilde r}_{1}={\tilde r}_{2}=\cdots={\tilde r}_{m}\in W$.

Finally, by the LaSalle's Theorem stated in Proposition 1, any  $r(t)$ starting from $C_\varepsilon$ converges to \(M\). From the properties of \(M\) analyzed above, it follows that $r(t) $  exhibits complete phase synchronization, and each \(r_{i}(t)\) approaches $W$ as $t\rightarrow +\infty$.   If the initial states satisfy $p^{\mathrm{T}}r_{i}(0)>0$ for some $p\in W$ and  every $i=1,2,\ldots,m$
, then we can choose sufficiently small $\varepsilon>0$ such that $(r(0),p)\in C_{\varepsilon}$. Therefore, the proof is complete.
\\ \ \par
Combining Theorem \ref{thm1} and Theorem \ref{thm2},  we get the necessary and sufficient condition for  complete phase synchronization: $\dim (W)>0$.

\begin{corollary}
        \label{cor4}
                Consider the high-dimensional Kuramoto model with nonidentical
                oscillators described by (\ref{eq001}). Assume that each $\Omega_i$ is skew-symmetric,   and the interconnection graph $\mathcal{G}$ is a strongly connected digraph. Then there is a unit semisphere $S$ such that all the solutions starting from $S$  exhibits complete phase synchronization if and only if $\dim (W)>0$, where $W$ is shown in (\ref{eqn006}).
\end{corollary}

\begin{remark}
        \label{rem2}
                If the nonidentical
high-dimensional Kuramoto oscillators are  proportional, i.e. $\Omega_i=l_i\Omega$ for each $i=1,2,\ldots,m$, then  $\dim (W)>0$ if and only if $\det (\Omega)=0$. Therefore, Corollary \ref{cor4} is a more general result than  the Theorem 3 of \cite{zhang2021synchronization}.
If $n$ is odd, then the skew-symmetric matrix $\Omega$ always satisfies $\det (\Omega)=0$.
\end{remark}
\begin{remark}
        \label{rem4}
In \cite{Markdahl2021CP}, Markdahl et al. proposed the condition that there exists a nonzero vector $v$ such that
$$
\mathrm{span}\{v\}=\bigcap_{i=1}^m\ker\Omega_i,
$$
which obviously implies that $W\neq \{0\}$. Our condition on the frequency matrices $\Omega_i$ is much weaker than the condition in \cite{Markdahl2021CP}. Moreover, compared with \cite{Markdahl2021CP}, in this paper the topology condition is weakened to strongly connected digraphs, and the constraint on the dimension $n$ is removed.
\end{remark}

\section{Simulations}

\ \ \ \ \ \  In this section, we give some simulations to illustrate the obtained theoretical results.
Consider the high-dimensional Kuramoto model (\ref{eq001}) composed by five oscillators connected by a
 directed cycle as follows:
 \begin{eqnarray}
\label{eq401}
\dot{r}_{1}&=&\Omega_{1}r_{1}+k(r_{5}-(r_{1}^{\mathrm{T}}r_{5})r_{1}),\nonumber  \\
\dot{r}_{2}&=&\Omega_{2}r_{2}+k(r_{1}-(r_{2}^{\mathrm{T}}r_{1})r_{2}),\nonumber \\
\dot{r}_{3}&=&\Omega_{3}r_{3}+k(r_{2}-(r_{3}^{\mathrm{T}}r_{2})r_{3}), \nonumber  \\
\dot{r}_{4}&=& \Omega_{4}r_{4}+k(r_{3}-(r_{4}^{\mathrm{T}}r_{3})r_{4}), \nonumber \\ \dot{r}_{5}&=& \Omega_{5}r_{5}+k(r_{4}-(r_{5}^{\mathrm{T}}r_{4})r_{5}). \nonumber
\end{eqnarray}

\begin{figure}[htb]
\label{fig2}
\centering
\includegraphics[height=7cm,width=9cm]{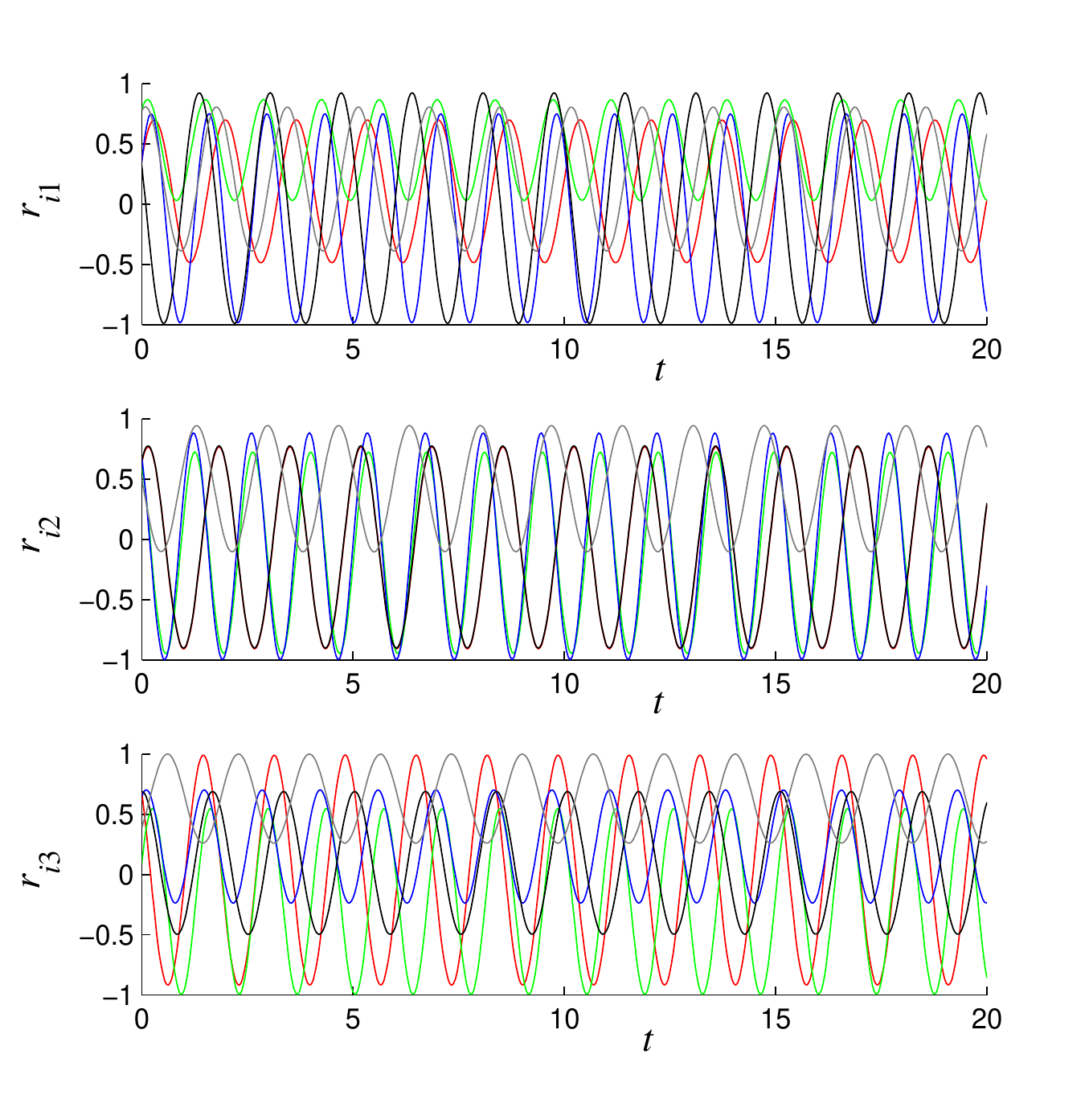}
\caption{{The time response curves of the states when $k=0$.}}
\label{fig402}
\end{figure}

\begin{figure}[htb]
\centering
\includegraphics[height=8cm,width=8.5cm]{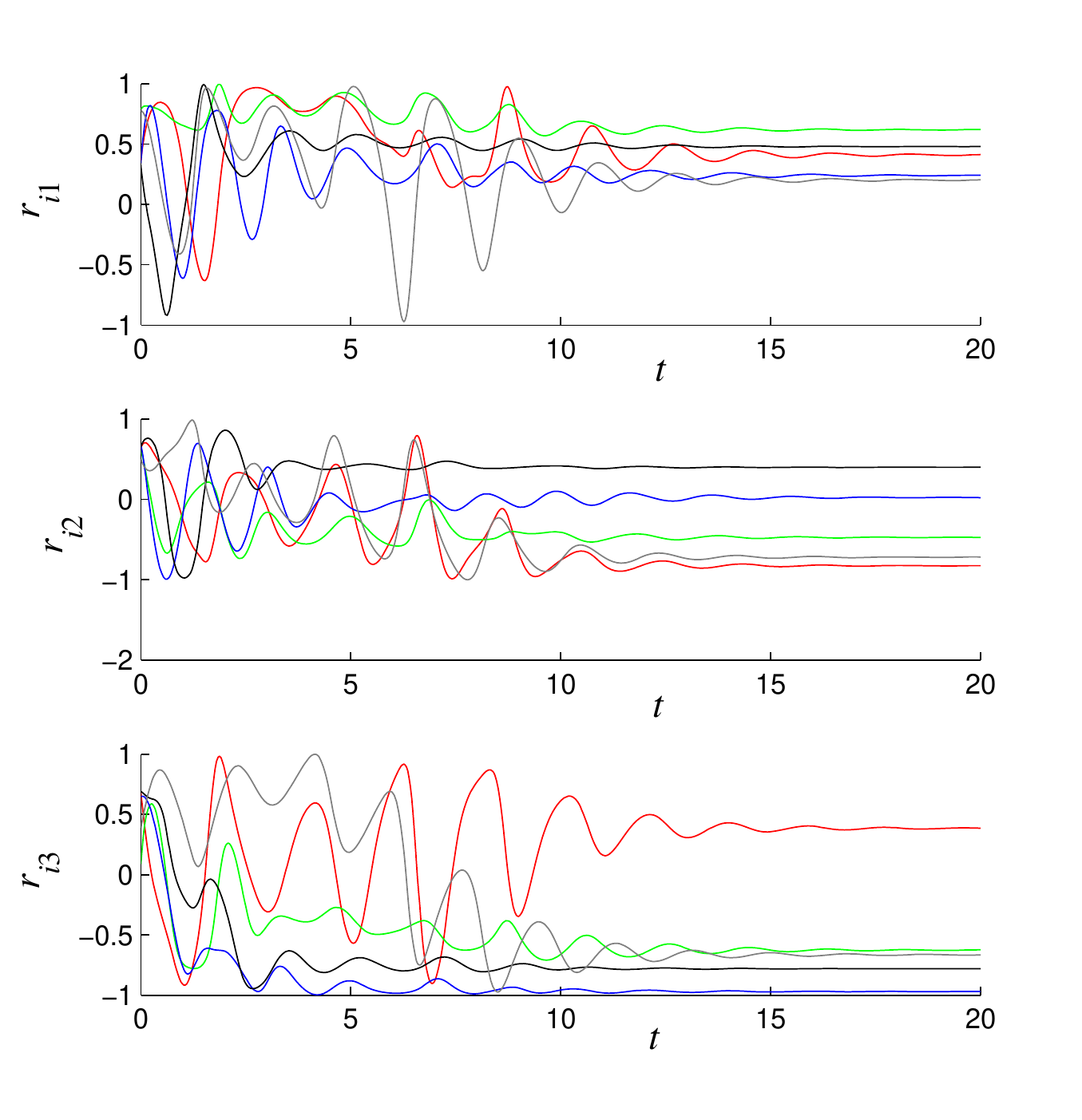}
\caption{{The time response curves of the states when $k=2$.}}
\label{fig404}
\end{figure}

\begin{figure}[htb]
\centering
\includegraphics[height=8cm,width=8cm]{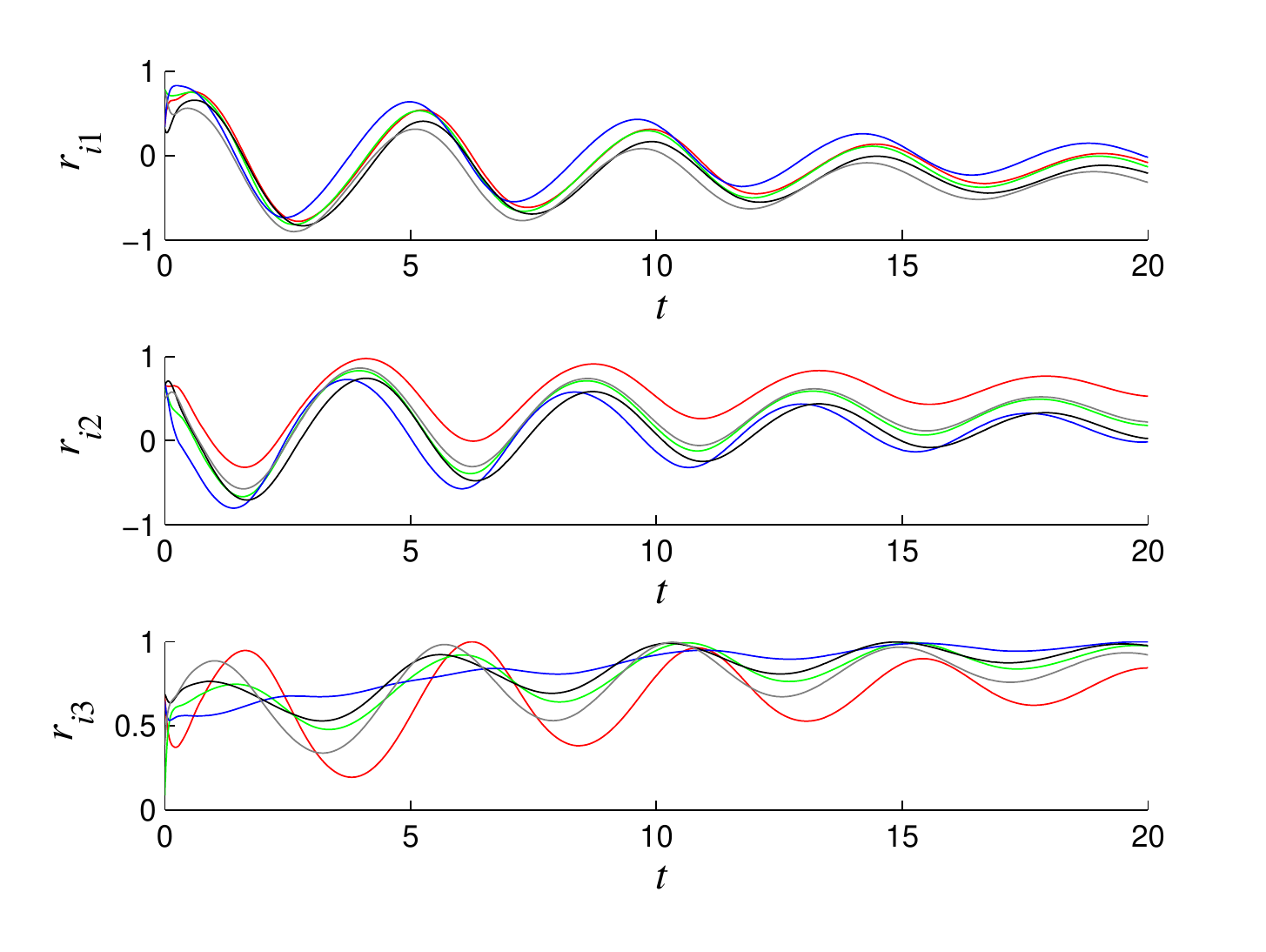}
\caption{{The time response curves of the states when $k=10$.}}
\label{fig404}
\end{figure}

First, we consider the model in 3-dimensional linear space, i.e. the case of $n=3$. Let
$$
\Omega_1\!\!=\!\!\left[\!\!\begin{array}{ccc} 0  &   1  &  2  \\
 -1&   0  &  3  \\
 -2&  -3  &  0
\end{array}\!\!\right],\
\Omega_2=\!\!\left[\!\!\begin{array}{ccc} 0   &  2 &   -1 \\                 -2  &  0 &   -4 \\
1   &  4 &    0
\end{array}\!\!\right],\
\Omega_3=\!\!\left[\!\!\begin{array}{ccc} 0   &  4 &   1 \\
-4  &  0 &   -2\\
-1  &  2 &   0
\end{array}\!\!\right],\
$$
$$
\Omega_4=\!\!\left[\!\!\begin{array}{ccc}   0  &  -3  &  -2 \\
3  &   0  &   1 \\
2  &  -1  &   0
\end{array}\!\!\right],\
\Omega_5=\!\!\left[\!\!\begin{array}{ccc}   0   &  3 &   -2\\
-3  &  0 &   1 \\
 2  &  -1&   0
\end{array}\!\!\right].
$$
It is easy to check that $$\mathrm{rank}\left[\begin{array}{c}
\Omega_2-\Omega_1\\
\Omega_3-\Omega_1      \end{array}\right]=3,$$ which implies that $$\ker(\Omega_2-\Omega_1)\bigcap \ker(\Omega_3-\Omega_1)=\{0\}.
$$ Therefore, we have $W=\{0\}$. By Theorem 1, we conclude that the complete phase synchronization cannot be achieved. We choose different $k$ for the simulations. In Figure 1-Figure 3, we draw the time response curves of the states when $k=0$, $k=2$ and $k=10$, respectively. The simulations shows the synchronization is not achieved even if $k$ is large, which validates our theoretical result of Theorem 1. However, the simulations show that the larger the value of $k$, the smaller the synchronization error. This is just the meaning of the so called practical synchronization.

\begin{figure}[htb]
\label{fig2}
\centering
\includegraphics[height=7cm,width=9cm]{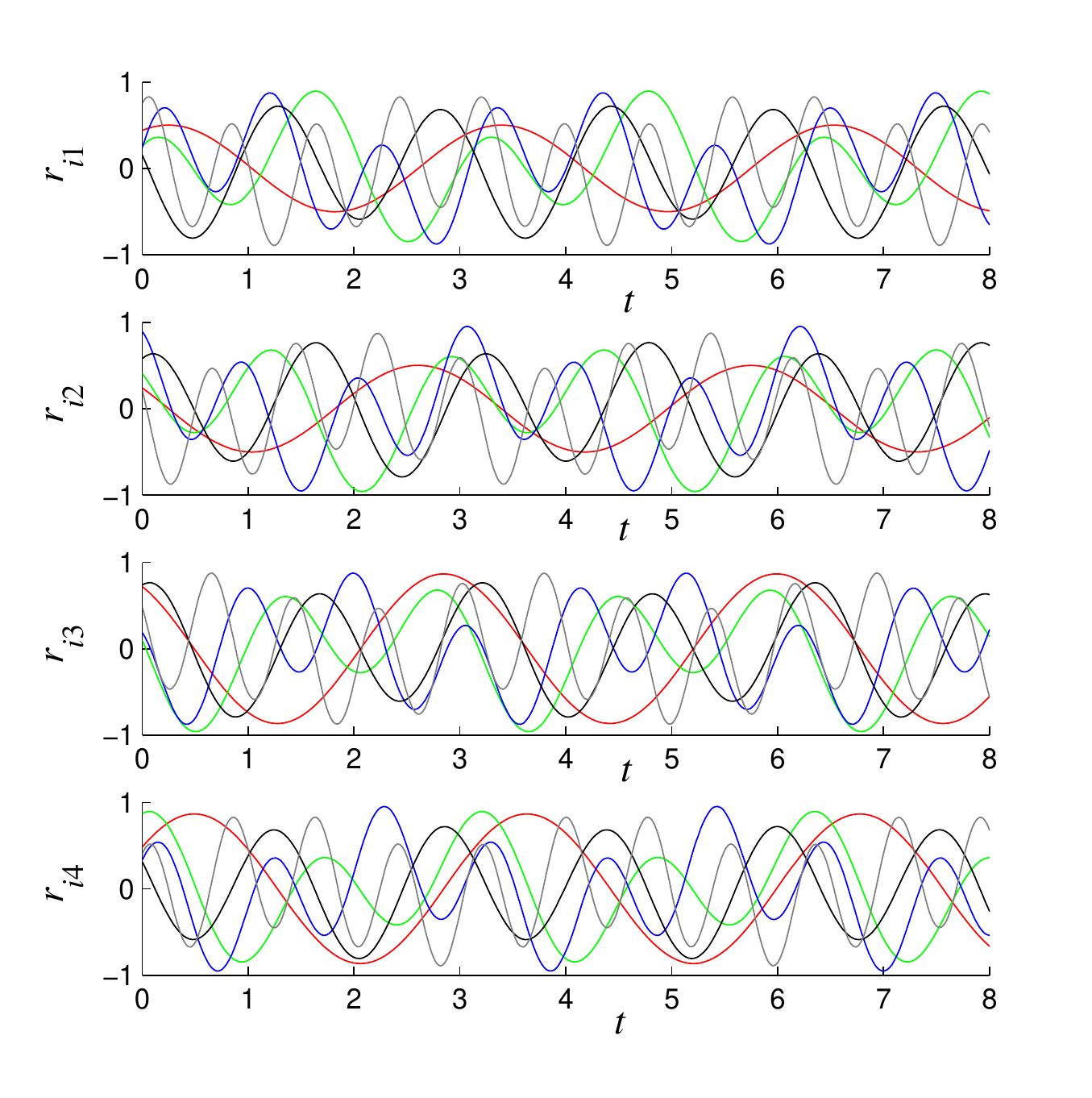}
\caption{{The time response curves of the states when $k=0$.}}
\label{fig402}
\end{figure}

\begin{figure}[htb]
\centering
\includegraphics[height=8cm,width=8.5cm]{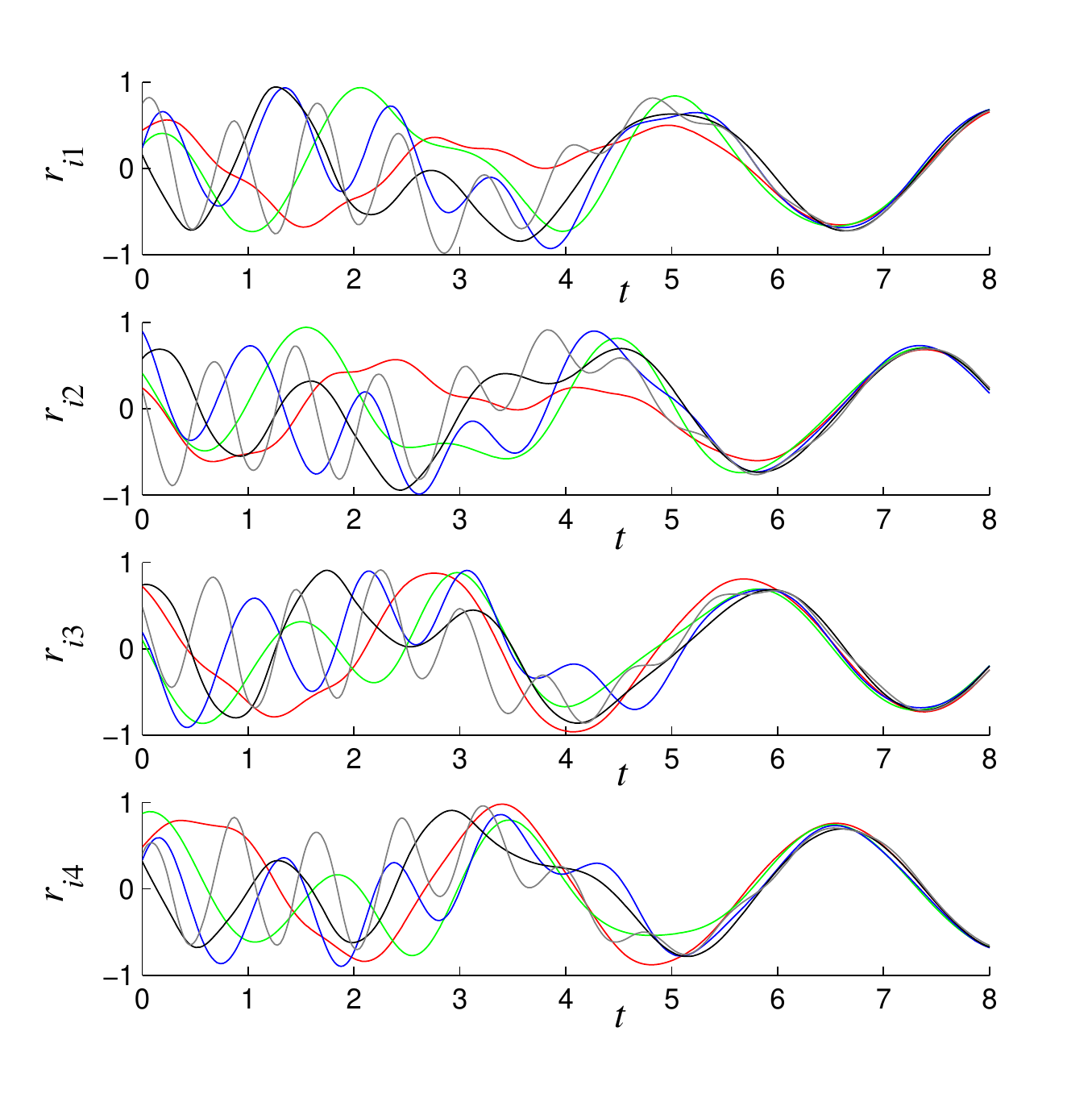}
\caption{{The time response curves of the states when $k=0.9$.}}
\label{fig404}
\end{figure}

\begin{figure}[htb]
\centering
\includegraphics[height=8cm,width=8cm]{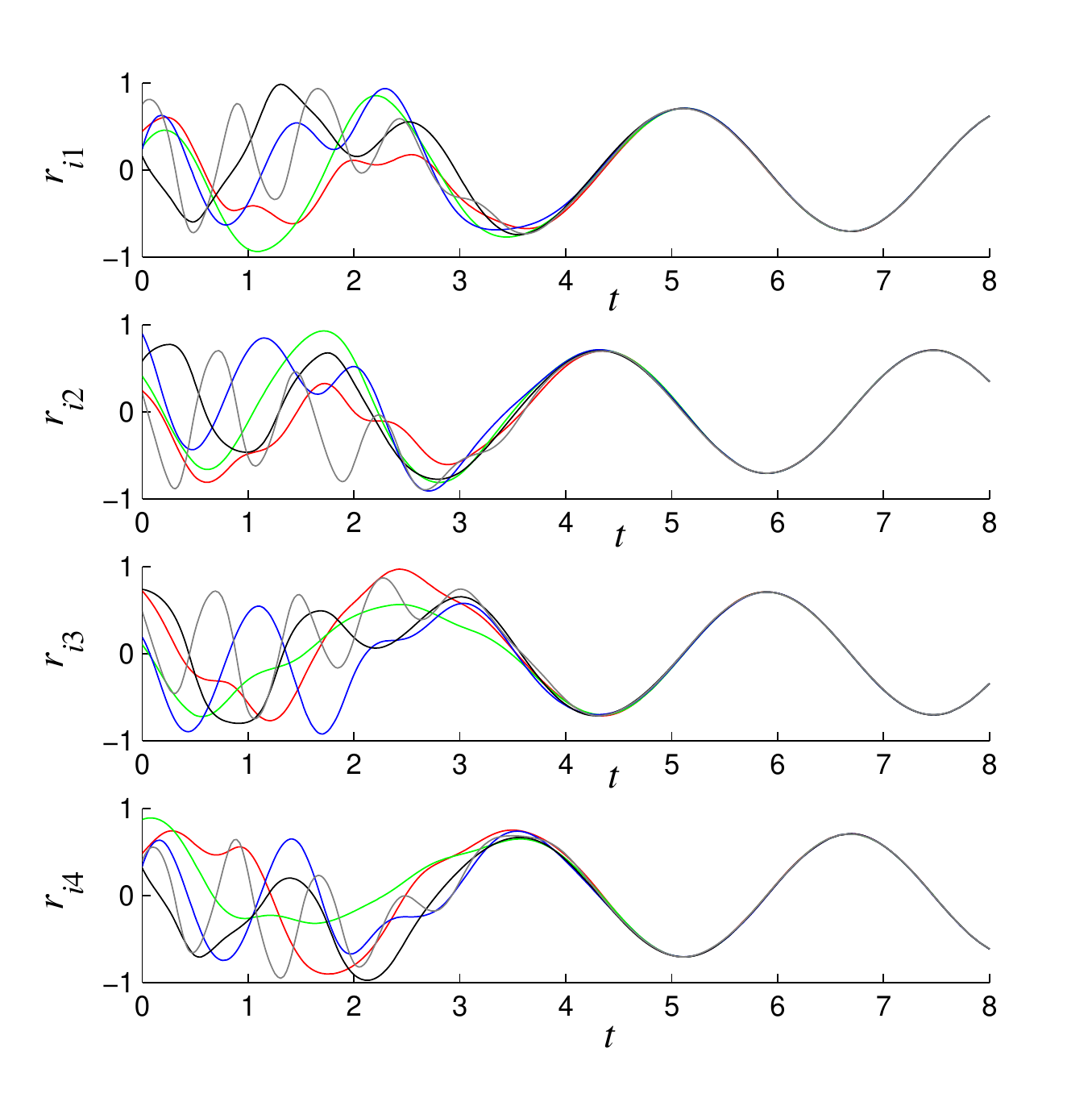}
\caption{{The time response curves of the states when $k=1.8$.}}
\label{fig404}
\end{figure}

Next, we consider the model in a 4-dimensional linear space.
Let
$$
\Omega_1=\!\!\left[\!\!\begin{array}{cccc}  0 & 2 &    0   &  0\\
-2& 0 &    0   &  0\\
 0& 0 &    0   & -2 \\
 0& 0 &    2   &  0
\end{array}\!\!\right],\
\Omega_2=\!\!\left[\!\!\begin{array}{cccc} 0&      3  &    1   &   0\\
-3&      0  &    0  &   -1\\
-1&      0  &    0  &   -3\\
 0&      1 &     3  &    0
\end{array}\!\!\right],\
\Omega_3=\!\!\left[\!\!\begin{array}{cccc} 0 &     4  &    2  &    0\\
-4&      0 &     0   &  -2\\
-2&      0&      0&     -4\\
 0&      2 &     4 &     0
\end{array}\!\!\right],
$$
$$
\Omega_4=\!\!\left[\!\!\begin{array}{cccc} 0&     -1 &    -3  &    0\\
  1&      0 &     0  &    3 \\ 3 &     0  &    0 &     1\\     0 &    -3  &   -1 &     0
\end{array}\!\!\right],\
\Omega_5=\!\!\left[\!\!\begin{array}{cccc} 0 &     5  &    3  &    0\\
 -5&      0&      0 &    -3\\
 -3&      0&      0  &   -5\\
  0&      3&      5  &    0
\end{array}\!\!\right].
$$
Since $\Omega_1^2=-4I_4$, from Lemma 5, it follows that
$$
W=\ker\left[\!\!\begin{array}{c}
\Omega_2-\Omega_1\\
\Omega_3-\Omega_1\\
\Omega_4-\Omega_1\\
\Omega_5-\Omega_1\\
(\Omega_2-\Omega_1)\Omega_1\\
(\Omega_3-\Omega_1)\Omega_1\\
(\Omega_4-\Omega_1)\Omega_1\\
(\Omega_5-\Omega_1)\Omega_1\\
\end{array}\!\!\right]=\mathrm{span}
\left[\!\!\begin{array}{cc}
0&1\\
1&0\\
-1&0\\
0&-1
\end{array}\!\!\right].
$$
Thus, by Theorem 2, the complete phase synchronization can be achieved. Under the initial conditions of Theorem 2, the states converge to $W$. In Figure 4-Figure 6,
we draw the time response curves of the states when $k=0$, $k=0.9$ and $k=1.8$, respectively. The simulations shows the synchronization is  achieved when $k>0$. From the limit state of each oscillator, we see that the first component plus the fourth one tends to zero, and the second component plus the third one also tends to zero, which means that each state converges to $W$. Meanwhile, the simulations show that the larger the value of $k$, the faster the synchronization.

Finally, we consider the model in 5-dimensional linear space, i.e. the case of $n=5$. Let
$$
\Omega_1\!\!=\!\!\left[\!\!\begin{array}{ccccc} 0   & 10 & -2  &  0 &  0\\
-10 & 0  &  4  &  0 &  0\\
2   & -4 &  0  &  0 &  0 \\
0   & 0  &  0  &  0 &  2 \\
0   & 0  &  0  & -2 &  0  \end{array}\!\!\right],\
\Omega_2\!\!=\!\!\left[\!\!\begin{array}{ccccc} 0& 1 & -7 &0& 0 \\                                      -1&0 & 1  &0& 0 \\
7 &-1& 0  &0& 0 \\
0 & 0& 0  &0& 2 \\
0 & 0& 0  & -2& 0
\end{array}\!\!\right],\
\Omega_3\!\!=\!\!\left[\!\!\begin{array}{ccccc} 0 & 3& 5& 0 & 0 \\
-3& 0& 2& 0 & 0\\
-5&-2& 0& 0 & 0 \\
0 & 0& 0& 0 & 2\\
0 & 0& 0&-2 & 0
\end{array}\!\!\right],
$$
$$
\Omega_4\!\!=\!\!\left[\!\!\begin{array}{ccccc} 0 & -4 & -1 & 0 & 0\\
4 & 0  & -1 & 0 & 0\\
1 & 1  &  0 & 0 & 0 \\
0 & 0  &  0 & 0 & 2\\
0 & 0  &  0 & -2& 0
\end{array}\!\!\right],\
\Omega_5\!\!=\!\!\left[\!\!\begin{array}{ccccc}
0 & 6 &  8 & 0 & 0\\
-6& 0 &  -3 & 0 & 0\\
-8& 3&  0 & 0 & 0\\
0 & 0 &  0 & 0 & -2\\
0 & 0 &  0 & 2 & 0\end{array}\!\!\right].
$$
It is easy to see that
\begin{equation}
\label{eqn017}
W=\mathrm{span}\{[0,0,0,1,0]^\mathrm{T},\ [0,0,0,0,1]^\mathrm{T}\}.
\end{equation}
\begin{figure}[htb]
\label{fig2}
\centering
\includegraphics[height=7cm,width=9cm]{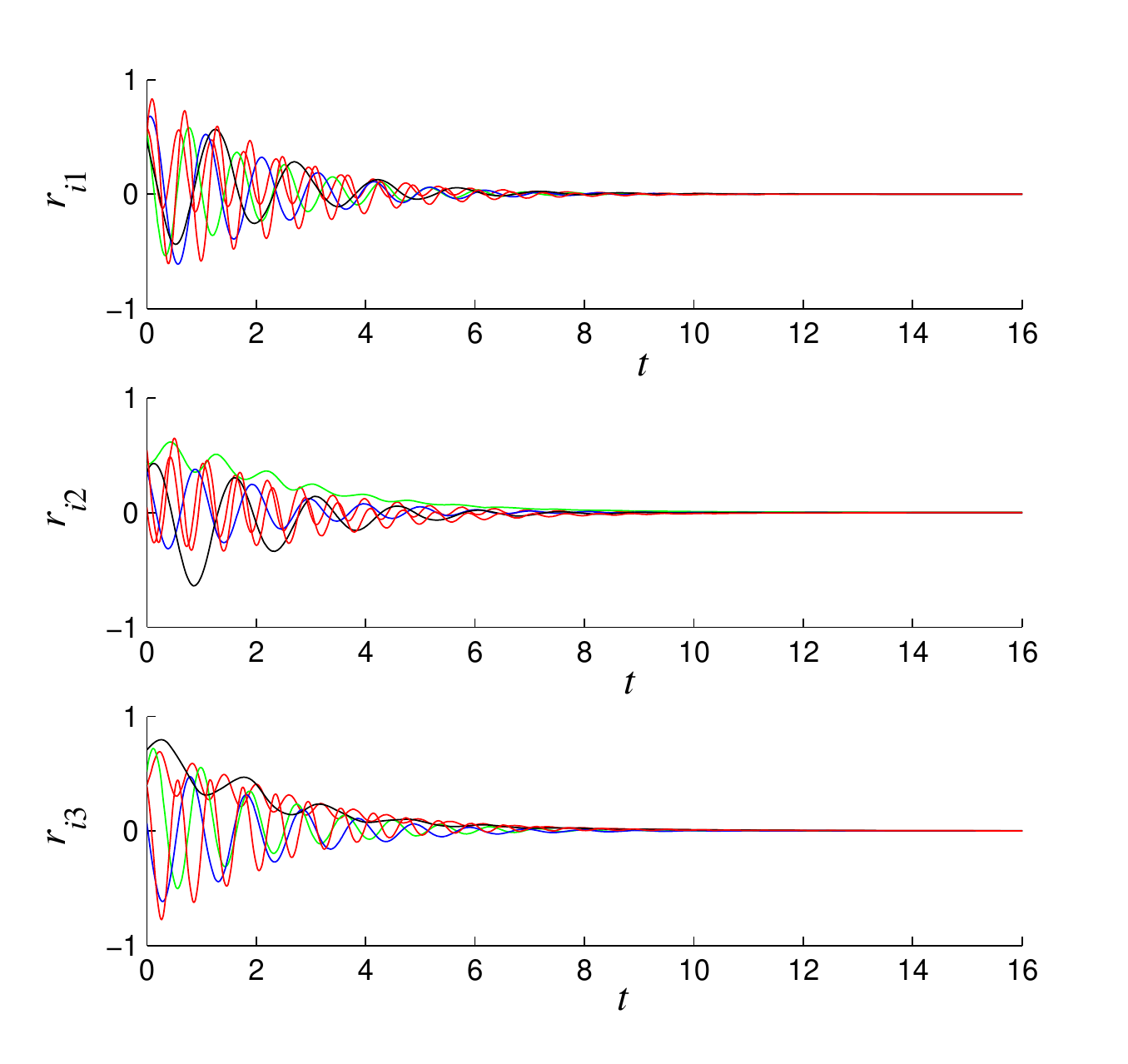}
\caption{{The time response curves of the first three components of the states.}}
\label{fig402}
\end{figure}

\begin{figure}[htb]
\centering
\includegraphics[height=8cm,width=8.5cm]{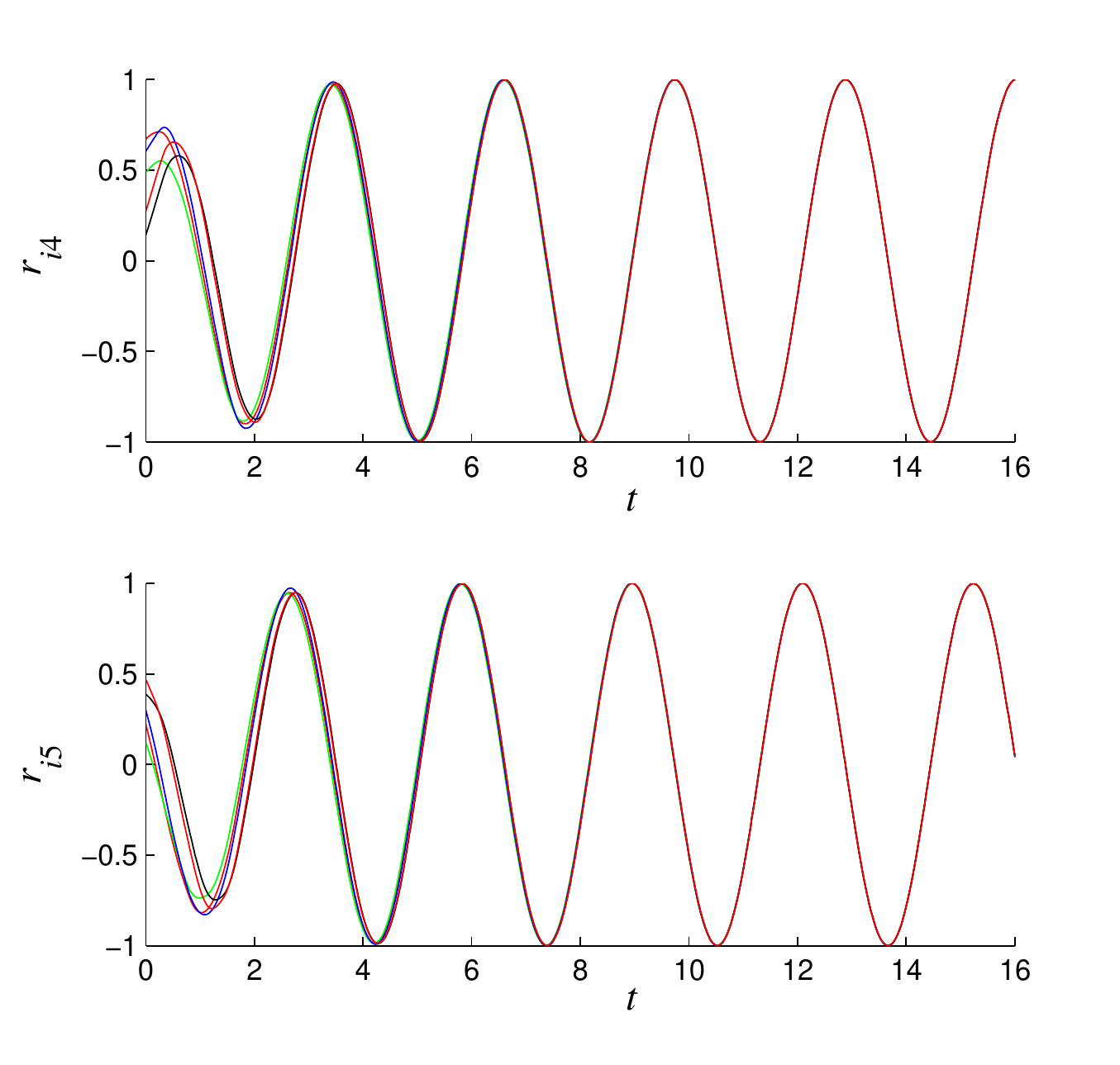}
\caption{{The time response curves of the fouth and the fifth components of the states.}}
\label{fig404}
\end{figure}

\begin{figure}[htb]
\centering
\includegraphics[height=8cm,width=8cm]{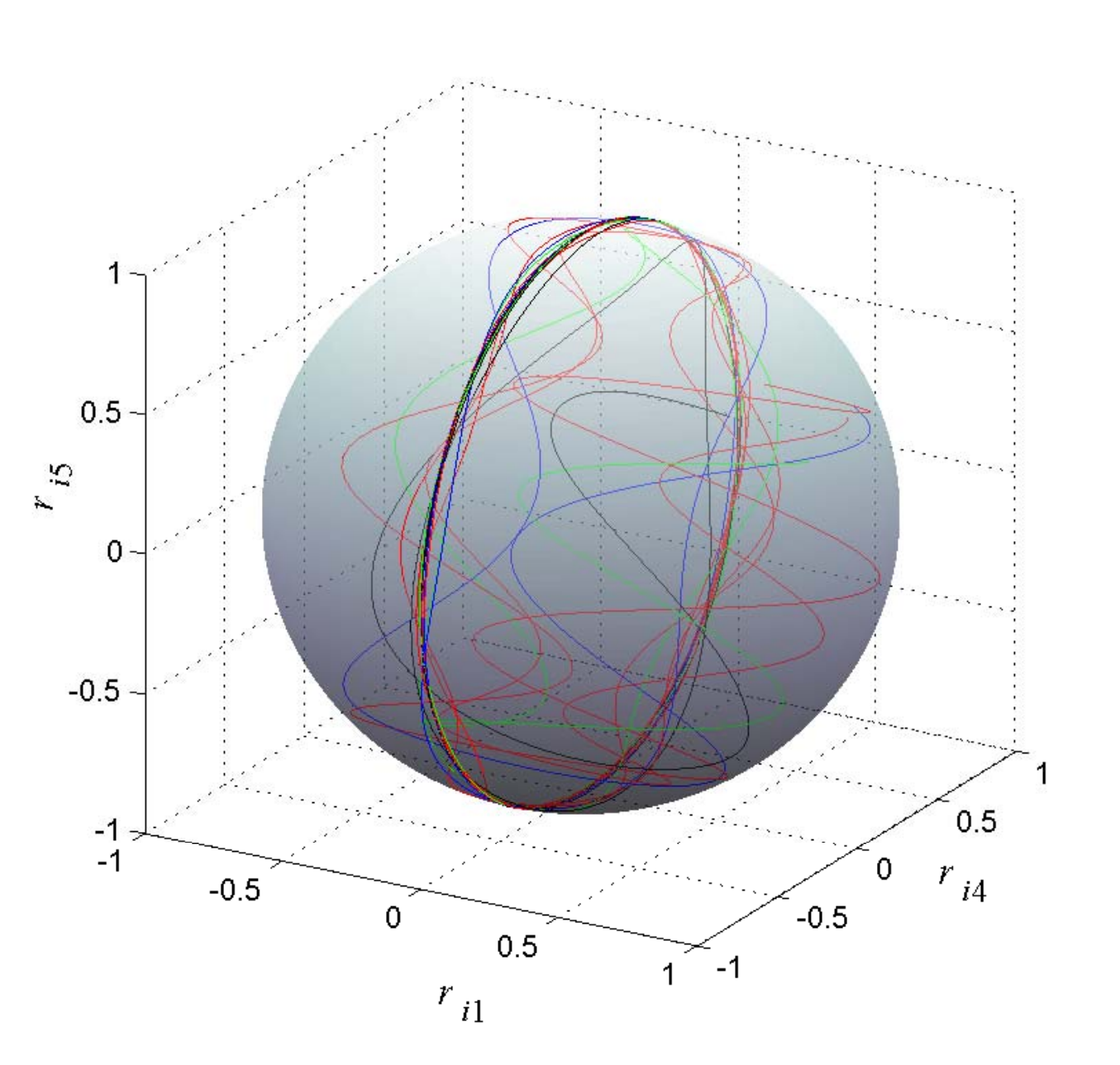}
\caption{{The the projection of the trajectories of the oscillators.}}
\label{fig404}
\end{figure}
 By Theorem 2, the complete phase synchronization can be achieved under the required initial condition, and the state of each oscillator converges to $W$. In Figure 7, it is displayed that the first three components of the oscillators' states tend to zero, which implies that $r_i(t)\rightarrow W$ ($t\rightarrow +\infty$) by  (\ref{eqn017}). Figure 8 shows that the complete phase synchronization is achieved.
Figure 9 shows the projection of the trajectories of the oscillators into the  three-dimensional space composed of the first, the fourth and the fifth components of the states.
The simulations validate the obtained theoretical results.

\section{Conclusion}
This paper has revealed that, unlike the original nonidentical Kuramoto model, a large class of high-dimensional nonidentical Kuramoto model can achieve complete phase synchronization. Under the topology of strongly connected digraphs, a necessary and sufficient condition for complete phase synchronization has been proposed. In our future work, we will focus on weakening the graph condition and strengthening the synchronization results.

\ \\ \  \\
{\bf References}
\ \\ \

\end{document}